\documentclass[11pt]{amsart}
\usepackage{amsmath}
\usepackage{amssymb}
\usepackage{amsthm}
\usepackage{amsrefs}
\usepackage{dsfont}
\usepackage{mathrsfs}
\usepackage{stmaryrd}
\usepackage{enumerate}
\usepackage[margin=1.1in]{geometry}
\usepackage[all]{xy}
\usepackage{verbatim}  
\usepackage{tikz, tikz-cd}
\usepackage{graphicx}
\graphicspath{ {./images_outputs/} }

\DeclareMathOperator{\Hom}{Hom}
\DeclareMathOperator{\tr}{tr}

\DeclareMathOperator{\Gal}{Gal}

\DeclareMathOperator{\Aut}{Aut}

\DeclareMathOperator{\Cay}{Cay}

\newcommand{\bC}{\mathbf{C}}

\newcommand{\bF}{\mathbf{F}}

\newcommand{\bQ}{\mathbf{Q}}

\newcommand{\bZ}{\mathbf{Z}}

\renewcommand{\phi}{\varphi}

\providecommand{\abs}[1]{\lvert #1 \rvert}

\providecommand{\oline}[1]{\overline{ #1 }}

\providecommand{\arr}{\longrightarrow}

\providecommand{\cn}{\colon}

\providecommand{\bFq}{\bF_q}
\providecommand{\bFsquares}{\bF_q^\boxtimes}
\providecommand{\bFsq}{\bFsquares}
\providecommand{\bFsp}{\bF_p^\boxtimes}
\providecommand{\bFplus}{\bF_q^+}
\providecommand{\bFtimes}{\bF_q^\times}

\providecommand{\resM}{M^{\boxtimes}}

\newtheorem{theorem}[equation]{Theorem}
\newtheorem{proposition}[equation]{Proposition}
\newtheorem{lemma}[equation]{Lemma}

\theoremstyle{definition}

\newtheorem{remark}[equation]{Remark}

\newtheorem{example}[equation]{Example}

\title[Cyclic covers of Paley graphs]{On the eigenvalues of cyclic covers of Paley graphs}
\author[Natalie Dinin]{Natalie Dinin}
\author[John A. Lind]{John A. Lind}
\address{California State University, Chico}
\email{jlind@csuchico.edu}

\begin{document}

\maketitle

\begin{abstract}
We study covering graphs of the Paley graph associated to a finite field of characteristic $p$ in the case where the covering transformation group is cyclic of prime order distinct from $p$.  When the field has $q = p$ elements, we show that the eigenvalues of the adjacency matrix determine the graph isomorphism class among translation invariant covers.  When $q = p^r > p$, we construct examples of cospectral covering graphs that are not isomorphic as graphs.
\end{abstract}

\section{Introduction}

The Paley graph $X(\bF_q)$ associated to a field $\bF_q$ with $q \equiv 1 \pmod{4}$ elements 
has vertex set $\bF_q$ and an edge between $x$ and $y$ if and only if $y - x \in \bFsq$, the set of quadratic residues in $\bF_q$.  The Paley graphs are strongly regular and self-complementary, and have been studied for their symmetry and embedding properties \citelist{\cite{Erdos_Renyi} \cite{Sachs} \cite{Bose} \cite{Jones:Paley}}.

In this paper, we study certain cyclic covering graphs of Paley graphs. 
We calculate the eigenvalues of their adjacency matrices and ask when covers with equivalent adjacency spectra are isomorphic as graphs.  In other words, in the spirit of Mark Kac \cite{Kac}, we ask: can one hear the shape of a cyclic cover of a Paley graph?  The answer is affirmative when $q = p$ is prime, 
but there are counterexamples when $q = p^r$ with $r \geq 2$.

In order to make a precise statement of the results, let us first describe the cyclic covering graphs under consideration.  Given a group $G$ and a connected graph $X$, a $G$-cover of $X$ is a connected graph $Y$ equipped with a free action of $G$ and an identification $Y/G \cong X$ of the quotient graph with $X$.  Such a covering graph is encoded by a function $\alpha$ that takes values in $G$ and is defined on the set of arcs in the base graph, meaning edges in $X$ with chosen orientation.  
For each closed path $e_1, \dotsc, e_n$ in $X$, the path-lifting monodromy action on the fibers of the covering graph $Y$ is given by the element $\alpha(e_1) \dotsm \alpha(e_n)$ of $G$.  See, for example, the book by Sunada \cite{Sunada} for background on the theory of covering graphs.

We are interested in the case where $X = X(\bF_q)$ is the Paley graph associated to a finite field $\bF_q$ of characteristic $p$, $G = \bZ/\ell$ is a cyclic group of prime order $\ell \neq p$, and $\alpha$ is invariant under the translation action by the additive group $\bFq^+$ on the arcs of $X(\bFq)$.  Under these hypotheses, the monodromy action is determined by a function $\bFsq \to \bZ/\ell$, which in abuse of notation we also denote by $\alpha$, that is not identically zero and satisfies $\alpha(-s) = - \alpha(s)$.  In this situation, we say that the associated covering graph $X^{\alpha}$ is a translation invariant $\bZ/\ell$-cover of $X(\bF_q)$.  Our main result says that when $q = p$ is prime, the adjacency eigenvalues of $X^{\alpha}$ determine its graph isomorphism class among such covers.

\begin{theorem}\label{thm:intro}
Suppose that $p$ is a prime number, and that $X^{\alpha}$ and $X^{\beta}$ are translation invariant $\bZ/\ell$-covers of the Paley graph $X(\bF_p)$ for a prime $\ell \neq p$.  If the adjacency spectra of $X^\alpha$ and $X^\beta$ are equal as multisets, then there is an isomorphism of graphs $X^{\alpha} \cong X^{\beta}$.  
\end{theorem}

\noindent For regular graphs, adjacency eigenvalues determine the eigenvalues of the graph Laplacian, and vice versa, so the same result holds for the Laplacian spectra of $X^{\alpha}$ and $X^{\beta}$.   

To prove the theorem, we first record basic facts about the Paley graphs in \S\ref{sec:paley_setup}, and then use character theory to show that the eigenvalues of the adjacency matrix of a covering graph $X^{\alpha}$ are given by the exponential sums
\begin{equation*}
\theta^{\alpha}_{a, k} = \sum_{s \in \bFsquares} \zeta_p^{\tr(as)}\zeta_{\ell}^{k\alpha(s)} \qquad \text{for $a \in \bFq, k \in \bZ/\ell$,}
\end{equation*}
where $\zeta_p = e^{2 \pi i /p}$ and $\zeta_{\ell} = e^{2 \pi i /\ell}$ are fixed choices of $p$-th and $\ell$-th roots of unity and $\tr \colon \bF_q^{+} \to \bF_p^{+}$ is the trace homomorphism. When $k = 0$, the eigenvalues $\theta^{\alpha}_{a, 0}$ are the adjacency eigenvalues of the Paley graph $X(\bF_q)$, and can be evaluated explicitly using the quadratic Gauss sum.  The eigenvalue calculations and the proof of Theorem \ref{thm:intro} take up \S\ref{sec:covers}.

When $q = p^r > p$, there are counterexamples to the statement in Theorem \ref{thm:intro}.  In \S\ref{sec:counterexamples}, we construct cospectral non-isomorphic covers for $q = 25$, $\ell > 2$, and discuss a general method for establishing the cospectrality of such examples.  In order to detect when covers are not isomorphic as graphs, we prove the following rigidity result in \S\ref{sec:CI_property}.

\begin{theorem}\label{thm:rigidity_intro}
Assume that $r \leq 4$ or $\ell > (q - 1)/2$.  If the translation invariant $\bZ/\ell$-covers $X^\alpha$ and $X^\beta$ of $X(\bFq)$ are isomorphic as graphs, then $X^{n \alpha}$ and $X^\beta$ are isomorphic as $\bZ/\ell$-covering graphs for some $n \in (\bZ/\ell)^\times$.
\end{theorem}
\noindent The proof amounts to establishing a version of the Cayley isomorphism property \cite{Babai:CI} for certain subsets of the group $\bF_q^+ \times \bZ/\ell$. A theorem of Carlitz \cite{Carlitz} implies that automorphisms of the Paley graph take a specific algebraic form, 
and a similar description can be deduced for isomorphisms of $\bZ/\ell$-covers.  Since the counterexamples do not admit an isomorphism of this form, Theorem \ref{thm:rigidity_intro} implies that they cannot be isomorphic as graphs.  


\medskip

{\bf Acknowledgements.}  We thank Thomas Mattman for introducing us to the Paley graphs.  We also thank Daniel Valli\`{e}res and Mooi Sprys for inspiring conversations about this work.

\numberwithin{equation}{section}
\section{The Paley graphs and their adjacency eigenvalues}\label{sec:paley_setup}

We summarize here the facts that we will need about Paley graphs.  This material can be found in many sources; for example, see \cite{Jones:Paley} for an overview and history of the Paley graphs and \cite{Nica} for the eigenvalue calculations.

Let $\bF_q$ be a field of characteristic $p$ with $q = p^r$ elements.  Write $\bFplus$ for the underlying additive group and $\bFtimes$ for the multiplicative group of nonzero elements.   We assume throughout that $q \equiv 1 \pmod{4}$.  Under this assumption, the group of quadratic residues $\bFsq = \{ x^2 \mid x \in \bFtimes\}$
contains the element $-1$, and so the $\bFq \times \bFq$-indexed matrix 
\begin{equation*}
A_{x, y} = \begin{cases} 1 \quad &\text{if $y - x \in \bFsq$,} \\ 0 \quad &\text{otherwise,} \end{cases}
\end{equation*}
is symmetric. By definition, the Paley graph $X(\bFq)$ has vertex set $\bFq$ and adjacency matrix $A$.  Equivalently, $X(\bFq) = \Cay(\bFplus, \bFsq)$ is the Cayley graph of the additive group $\bFplus$ with respect to the generating set $\bFsq \subset \bFplus$.  

As is true for any Cayley graph, the characters $\chi \in \bF_q^\vee = \Hom(\bF^+_q, \bC^\times)$ are eigenfunctions for the adjacency operator determined by $A$:
\[
(A\chi)(x) = \sum_{\substack{\text{edges $x \sim y$} \\ \text{in $X(\bFq)$}}} \chi(y) = \sum_{s \in \bF_q^\boxtimes} \chi(x + s) = \sum_{s \in \bF_q^\boxtimes} \chi(s) \cdot \chi(x).
\]
To calculate the adjacency eigenvalues explicitly, use the fact that all characters of $\bFplus$ are of the form
\[
\chi_a(x) = \zeta_p^{\tr(ax)} \qquad \text{for $a \in \bF_q$,}
\]
where we fix the choice $\zeta_p = e^{2\pi i/p}$ for a primitive $p$-th root of unity, and $\tr \colon \bF_q^+ \arr \bF_p^+$ is the additive trace homomorphism $\tr(x) = x + x^p + \dotsm + x^{p^{r - 1}}$.
Thus, the adjacency eigenvalues of the Paley graph $X(\bFq)$ are given by the character sums 
\begin{equation}\label{eq:theta_a}
\theta_a = \sum_{s \in \bF_q^\boxtimes} \chi_a(s) = \sum_{s \in \bF_q^\boxtimes} \zeta_p^{\tr(as)} \quad \text{for $a \in \bF_q$.}
\end{equation}   
When $a = 0$, we have the trivial character $\chi_0$ with adjacency eigenvalue $\theta_0 = (q - 1)/2$.  For non-trivial characters $\chi_a$, one relates the sum over squares in \eqref{eq:theta_a} to the Gauss sum $\sum \eta(x) \chi_a(x)$ for the quadratic character $\eta$ of $\bFtimes$.  The adjacency eigenvalues in this case are $(\pm \sqrt{q} - 1)/2$, each with multiplicity $(q - 1)/2$.  More explicitly, we have
\begin{equation}\label{eq:gauss_sum_calculation}
\theta_a = \begin{cases}
\displaystyle \sum_{s \in \bFsq} \zeta_p^{\tr(s)} = \frac{(-1)^{r - 1}\sqrt{q} - 1}{2} \quad &\text{when $a \in \bFsq$,} \\ 
\displaystyle \sum_{s \in \bFtimes \setminus \bFsq} \zeta_p^{\tr(s)} = \frac{(-1)^{r}\sqrt{q} - 1}{2} \quad &\text{when $a \in \bFtimes \setminus \bFsq$.}
\end{cases} 
\end{equation}
When $q = p$, this goes back to Gauss.  When $q = p^r > p$, a theorem of Davenport and Hasse \cite{Davenport_Hasse} reduces the calculation to the classical result of Gauss (see \S11.3, Thm. 1 in \cite{Ireland_Rosen}).


The group of automorphisms of the Paley graph $X(\bF_q)$ is very rigid.  The only possible automorphisms act on the vertex set $\bF_q$ by
\begin{equation}\label{eq:Paley_autos}
x \mapsto tx^\sigma + a  \quad \text{for some } t \in \bFsq, \sigma \in \Gal(\bF_q/\bF_p), a \in \bF_q.
\end{equation}
This result is usually attributed to Carlitz, since it follows directly from his theorem \cite{Carlitz} stating that any bijective function $f \cn \bFq \to \bFq$ satisfying $f(0) = 0$,  $f(1) = 1$, and 
\begin{equation*}
\eta(f(x) - f(y)) = \eta(x - y),
\end{equation*}
where $\eta \cn \bFtimes \to \{\pm 1\}$ is the quadratic character, must be of the form $f(x) = x^{p^i}$ for some $0 \leq i \leq r - 1$.  See \cite{Jones:Paley}*{\S9} and the references therein for generalizations and further commentary.

\section{Cyclic covers of Paley graphs and their adjacency eigenvalues}\label{sec:covers}

We now introduce the cyclic covers $X^\alpha$ of the Paley graph $X(\bFq)$ that are our primary objects of study.  Fix a prime number $\ell \neq p$ and let $\bZ/\ell$ be the cyclic group of integers under addition modulo $\ell$.  Suppose that $\alpha \cn \bF_q^{\boxtimes} \to \bZ/\ell$ is a function with $\alpha(s) \neq 0$ in $\bZ/\ell$ for some $s \in \bFsq$ and satisfying $\alpha(-s) = -\alpha(s)$ for all $s \in \bFsq$.  We call such a function $\alpha$ (which is generally not a homomorphism with respect to either addition or multiplication in $\bFq$) a voltage assignment.

Let $X^{\alpha}$ be the graph with vertex set $\bF_q \times \bZ/\ell$ and an edge between $(x, i)$ and $(y, j)$ if and only if $y - x \in \bFsq$ and $j - i = \alpha(y - x)$.  Then the projection to the first coordinate makes $X^\alpha \to X(\bF_q)$ a $\bZ/\ell$-covering graph of the Paley graph.  Our assumptions guarantee that the $\bZ/\ell$ action on the cover $X^\alpha$ may be promoted to an action of the group $\bF^+_q \times \bZ/\ell$ that is compatible, via the projection map, with the action of $\bFplus$ on $X(\bFq)$ by translation.  Every covering graph of $X(\bFq)$ with this property, which we call a translation invariant $\bZ/\ell$-cover, arises as $X^\alpha$ for some voltage assignment $\alpha$.

The quotient graph $X^\alpha/(\bF^+_q \times \bZ/\ell)$ is a bouquet with $(q - 1)/4$ edges.  Equivalently, the covering graph $X^\alpha$ is the Cayley graph of the group $\bF^+_q \times \bZ/\ell$ with respect to the generating set 
\begin{equation*}
S^\alpha = \{(s, \alpha(s)) \mid s \in \bFsquares \} \subset \bF^+_q \times \bZ/\ell.
\end{equation*}
Just as for the Paley graphs, the eigenvalues of the adjacency matrix $A^{\alpha}$ of $X^{\alpha}$ are studied using character theory.  
Each character of the group $\bF_q^+ \times \bZ/\ell$ may be written as a product $(\chi, \psi)(x, i) = \chi(x)\cdot\psi(i)$ of characters $\chi \in (\bF_q^+)^\vee$ and $\psi \in (\bZ/\ell)^\vee$, and is an eigenfunction of the adjacency operator: 
\begin{align}\label{eq:char_eigenfunction}
A^\alpha(\chi, \psi)(x, i) = \sum_{s \in \bFsq} (\chi, \psi)(x + s, i + \alpha(s)) = \bigl( \sum_{s \in \bFsq} \chi(s)\psi(\alpha(s) )\bigr) (\chi, \psi)(x, i).
\end{align}
To write formulas for the eigenvalues, use the description $\chi_a(x) = \zeta_p^{\tr(ax)}$ for $a \in \bF_q$ of the characters of $\bFplus$ from \S\ref{sec:paley_setup} and, in addition, fix the choice $\zeta_\ell = e^{2 \pi i /\ell}$ of primitive $\ell$-th root of unity, so that each character of $\bZ/\ell$ is of the form $\psi_k(i) = \zeta_\ell^{ki}$ for $k \in \bZ/\ell$.  The adjacency eigenvalues of the graph $X^\alpha$ are the exponential sums
\begin{equation}\label{eq:eigenvalue_formula}
\theta^{\alpha}_{a, k} = \sum_{s \in \bFsquares} \chi_a(s)\psi_k(\alpha(s)) = \sum_{s \in \bFsquares} \zeta_p^{\tr(as)}\zeta_{\ell}^{k\alpha(s)} \quad \text{for $a \in \bFq$, $k \in \bZ/\ell$.}
\end{equation}

The trivial character ($a = 0, k = 0$) gives the maximal eigenvalue $\theta^{\alpha}_{0, 0} = \abs{\bFsq} = (q - 1)/2$,
which is the valence of the regular graph $X^\alpha$. Our assumption that $X^\alpha$ is connected implies that the eigenvalue $(q - 1)/2$ has multiplicity one.  
For nonzero $a \in \bF_q$ and $k = 0$, it is clear that $\theta^{\alpha}_{a, 0} = \theta_a$ is the corresponding eigenvalue of the Paley graph, with explicit values given by \eqref{eq:gauss_sum_calculation}.
Thus, the multiset $\{\theta_{a, 0}\}_{a \in \bF_q}$ recovers the spectrum of the Paley graph $X(\bF_q)$.

The condition $\alpha(-s) = -\alpha(s)$ implies that $\theta^\alpha_{-a, -k} = \theta^\alpha_{a, k}$.  When $\ell$ is odd, this means that the eigenvalue multiplicities, besides that of the maximal eigenvalue $\theta_{0, 0} = (q - 1)/2$, are even.  Splitting the sum of roots of unity up over pairs $\pm s \in \bFsq$ gives a slightly different expression
\begin{equation*}
\theta^{\alpha}_{a, k} = \sum_{\{s, -s\} \subset \bFsquares} (\zeta_p^{\tr(as)}\zeta_{\ell}^{k\alpha(s)} + \zeta_p^{\tr(-as)}\zeta_{\ell}^{k\alpha(-s)}) = \sum_{\{s, -s\} \subset \bFsquares} 2 \mathrm{Re} ( \zeta_p^{\tr(as)}\zeta_{\ell}^{k\alpha(s)} )
\end{equation*}
which makes explicit the fact that the eigenvalues are real (as must hold for any undirected graph).  

Characters form an orthonormal basis $\{\chi_a\}$ for the space of complex functions on $\bFq$ with respect to the inner product
\begin{equation*}
\langle \phi_1, \phi_2 \rangle = \frac{1}{q} \sum_{a \in \bFq} \phi_1(a) \overline{\phi_2(a)}.
\end{equation*}
The eigenvalues of $X^\alpha$ then satisfy $\theta^{\alpha}_{a, k} = q \cdot \langle \chi_a , \overline{\psi_{k}^{\alpha}} \rangle$, where for each character $\psi \in (\bZ/\ell)^\vee$ we write $\psi^\alpha \cn \bFq \to \bC$ for the function 
\begin{equation*}
\psi^\alpha(x) = \begin{cases}
\psi(\alpha(x)) \quad &\text{if $x \in \bFsq$,} \\ 
0 \quad &\text{else.}
\end{cases}
\end{equation*}
Some readers might prefer to think of $\theta^{\alpha}_{-, k}$ as the discrete Fourier transform of the function $\psi^{\alpha}_k$ over the finite abelian group $\bFplus$.  The next lemma, a direct consequence of the orthonormality of the basis $\{\chi_a\}$, then amounts to taking the inverse Fourier transform.

\begin{lemma}\label{lem:fourier_inverse}  
Let $k \in (\bZ/\ell)^\times$.
If $\alpha$ and $\beta$ are voltage assignments satisfying $\theta^{\alpha}_{a, k} = \theta^{\beta}_{a, k}$ for all $a \in \bFq$, then $\alpha = \beta$.
\end{lemma}

\medskip

We need one more result before beginning the proof of the main theorem from the introduction.

\begin{proposition}\label{prop:iso_covers}
Every isomorphism of $\bZ/\ell$-covering graphs $f \cn X^\alpha \to X^\beta$ is of the form 
\begin{equation*}
f(x, i) = (t x^{\sigma} + a, i + k) 
\end{equation*}
for some $t \in \bFsquares, \sigma \in \Gal(\bF_q/\bF_p), a \in \bF_q$, $k \in \bZ/\ell$ satisfying $\alpha(s) = \beta(t s^\sigma)$ for all $s \in \bFsquares$.  Conversely, any function of the indicated form determines an isomorphism of $\bZ/\ell$-covering graphs. 
\end{proposition}
\begin{proof}
An isomorphism of $\bZ/\ell$-covering graphs respects the projections onto the base graph, and so fits into a commutative diagram 
\begin{equation*}
\begin{tikzcd}
 X^\alpha \arrow[r, "f"] \arrow[d, ""] & X^\beta \arrow[d, ""] \\ 
 X(\bF_q) \arrow[r, "\oline{f}"] & X(\bF_q)
\end{tikzcd}
\end{equation*}
with $\oline{f} \cn X(\bF_q) \to X(\bF_q)$ an automorphism of the Paley graph $X(\bF_q)$.  According to the description \eqref{eq:Paley_autos} afforded by Carlitz's theorem \cite{Carlitz}, we can find $t \in \bFsquares, \sigma \in \Gal(\bF_q/\bF_p),$ and $a \in \bF_q$ such that $\oline{f}(x) = tx^\sigma + a$.  Restricting to inputs with $i = 0$, we have 
\begin{equation*}
f(x, 0) = (\oline{f}(x), k(x)) = (tx^\sigma + a, k(x))
\end{equation*} 
for some function $k \cn \bF_q \to \bZ/\ell$.  Being a map of $\bZ/\ell$-covers, $f$ respects the translation action of $\bZ/\ell$, and so 
\begin{equation*}
f(x + s, \alpha(s)) = f(x + s, 0) + (0, \alpha(s)) = (t x^\sigma + ts^\sigma + a, k(x + s) + \alpha(s))
\end{equation*}
for any $x \in \bF_q$, $s \in \bFsq$.  The edge connecting $(x, 0)$ and $(x + s, \alpha(s))$ in $X^\alpha$ is mapped under $f$ to an edge in $X^\beta$ and so 
\begin{equation*}
\beta(t s^\sigma) = k(x + s) - k(x) + \alpha(s).
\end{equation*}
Taking the sum of these relations in $\bZ/\ell$ for $x = 0, s, \dotsc, (p - 1)s$ gives the equation 
\begin{equation*}
p(\beta(t s^\sigma) - \alpha(s)) = 0
\end{equation*}
from which we conclude that $\beta(t s^\sigma) = \alpha(s)$.  Thus, $k(x + s) = k(x)$ for all $s \in \bFsq$.  Since $\bFsq$ generates $\bFplus$, it follows that $k$ is a constant function $k(x) = k$ and so $f$ takes the form in the statement of the lemma.  A quick check shows that any function of this form defines an isomorphism of $\bZ/\ell$-covers.
\end{proof}

\begin{remark}
It is convenient to use the notation $\beta^{t}$ for the voltage assignment $s \mapsto \beta(t^{-1} s)$.  A consequence of the proposition is that the graphs $X^\alpha$ and $X^\beta$ are isomorphic as $\bZ/\ell$-covers if and only if $\alpha = \beta^{t} \circ \sigma$ for some $t \in \bFsquares, \sigma \in \Gal(\bF_q/\bF_p)$.  Given $n \in (\bZ/\ell)^{\times}$, the graph $X^{n\alpha}$ associated to the voltage assignment $s \mapsto n \cdot \alpha(s)$ is isomorphic to $X^{\alpha}$, via the isomorphism $X^{\alpha} \to X^{n\alpha}$ given by $(x, i) \mapsto (x, ni)$, but $X^{\alpha}$ and $X^{n\alpha}$ are not isomorphic as $\bZ/\ell$-covers unless $n = \pm 1$.

One might ask if every isomorphism of graphs $X^\alpha \to X^\beta$ is of the form given in the proposition.  This amounts to determining if, up to translation by elements of the group $\bF^+_q \times \bZ/\ell$, every isomorphism may be written $f(x, i) = (tx^\sigma, i)$.  As stated, the request is clearly false, because group automorphisms of $\bZ/\ell$ can also act on the second coordinate.  It is more reasonable to ask that $f(x, i) = (tx^\sigma, ni)$ for some $n \in (\bZ/\ell)^\times$.  In \S\ref{sec:CI_property}, we give an affirmative answer to this question under the additional assumption that $(q - 1)/2 < \ell$. 
\end{remark}


We now give an expanded statement and the proof of Theorem \ref{thm:intro} from the introduction.

\begin{theorem}
Suppose that $q = p$ is a prime number, and that $X^\alpha$ and $X^\beta$ are translation invariant $\bZ/\ell$-covers of the Paley graph $X(\bF_p)$.  If the adjacency spectra of $X^\alpha$ and $X^\beta$ coincide, meaning that $\{\theta^\alpha_{a, k}\} = \{\theta^\beta_{a, k}\}$ as multisets, then there is an isomorphism of $\bZ/\ell$-covering graphs $X^\alpha \cong X^{n \beta}$ for some $n \in (\bZ/\ell)^\times$, and so in particular $X^\alpha$ and $X^\beta$ are isomorphic as graphs. 
\end{theorem}
\begin{proof}
Assume that the adjacency spectra of $X^\alpha$ and $X^\beta$ agree.  Then $\theta^\alpha_{1, 1} = \theta^\beta_{t, n}$ for some $t \in \bF_p$, $n \in \bZ/\ell$.  Since the eigenvalues $\theta^\alpha_{a, 0} = \theta^\beta_{a, 0}$ agree for all $a \in \bF_p$, we may assume that $n \neq 0$.  The bulk of the proof is in reducing to the case $t \in \bFsp$.  

To this end, consider the given equation of eigenvalues
\begin{equation}\label{eq:eigenvalue_1_1_assumption}
\theta^\alpha_{1, 1} = \sum_{s \in \bFsp} \zeta_p^{s}\zeta_{\ell}^{\alpha(s)} = \sum_{s \in \bFsp} \zeta_p^{ts}\zeta_{\ell}^{n\beta(s)} = \theta^{\beta}_{t, n}
\end{equation}
and the effect of applying the field automorphism $\bQ(\zeta_p, \zeta_\ell) \to \bQ(\zeta_p, \zeta_\ell)$ induced by $\zeta_\ell \mapsto \zeta_\ell^k$ for $k \in (\bZ/\ell)^\times$.  Since $p$ is prime to $\ell$, there is no effect on $\zeta_p$ and we have
\begin{equation*}
\sum_{s \in \bFsp} \zeta_p^{s}\zeta_{\ell}^{k\alpha(s)} = \sum_{s \in \bFsp} \zeta_p^{ts}\zeta_{\ell}^{kn\beta(s)}.
\end{equation*}
Summing over $k$ and using the relation $\displaystyle\sum_{i \neq 0} \zeta_\ell^i = -1$ gives
\begin{equation*}
(\ell - 1)\sum_{\alpha(s) = 0} \zeta_p^s - \sum_{\alpha(s) \neq 0} \zeta_p^s = (\ell - 1)\sum_{n\beta(s) = 0} \zeta_p^{ts} - \sum_{n\beta(s) \neq 0} \zeta_p^{ts},
\end{equation*}
where we have split the terms up according to whether $\zeta_{\ell}^{k\alpha(s)}, \zeta_{\ell}^{kn \beta(s)} = 1$ or not.  It follows that
\begin{equation*}
\sum_{s \in \bFsp} \zeta_p^s \equiv \sum_{s \in \bFsp} \zeta_p^{ts} \quad \pmod{\ell \bZ[\zeta_p]}.
\end{equation*}
The term on the left is the Paley eigenvalue $(\sqrt{p} - 1)/2$.  The term on the right is either $(p - 1)/2$ or $(\pm \sqrt{p} - 1)/2$, according to whether $t = 0$, $t \in \bFsp$, or $t \in \bF_p^{\times} \setminus \bFsp$.  
In the last case, when $t \in \bF_p^{\times} \setminus \bFsp$, we have 
\begin{equation*}
\frac{\sqrt{p} - 1}{2} \equiv \frac{-\sqrt{p} - 1}{2} \pmod{\ell\bZ[\zeta_p]},
\end{equation*}
and so $\sqrt{p} \in \ell\bZ[\zeta_p]$. But then $\ell$ divides $\sqrt{p}$, and so the decompositions of $(p)$ and $(l)$ into prime ideals in $\bZ[\zeta_p]$ share a common factor, which is absurd.  In the case when $t = 0$, we have $(\sqrt{p} - 1)/2 \equiv (p - 1)/2$ mod $\ell\bZ[\zeta_p]$, and so $\ell$ divides $\sqrt{p}(\sqrt{p} - 1)/2$ in the ring $\bZ[\zeta_p]$.  It follows that $\ell$ must divide $(\sqrt{p} - 1)(\sqrt{p} + 1)/2 = (p - 1)/2$ in $\bZ$.  To finish off the case $t = 0$, return to the assumption \eqref{eq:eigenvalue_1_1_assumption} and apply the field automorphism $\bQ(\zeta_p, \zeta_\ell) \to \bQ(\zeta_p, \zeta_\ell)$ induced by $\zeta_p \mapsto \zeta_p^a$ for $a \in (\bZ/p)^\times$:
\begin{equation}\label{eq:eigenvalue_1_1_galois_act}
\sum_{s \in \bFsp} \zeta_p^{as}\zeta_{\ell}^{\alpha(s)} = \sum_{s \in \bFsp} \zeta_p^{ats}\zeta_{\ell}^{n\beta(s)} = \sum_{s \in \bFsp} \zeta_{\ell}^{n\beta(s)}.
\end{equation}
The resulting equation $\theta^\alpha_{a, 1} = \theta^\beta_{0, n}$ holds for $a = 1, \dotsc, p - 1$ and so the multiplicity of the eigenvalue $\theta^\beta_{0, n}$ in the adjacency spectrum for $X^\beta$ must be at least $p - 1$.  Since there are only $\ell$ eigenvalues of the form $\theta^\beta_{0, k}$ and $\ell \mid (p - 1)/2 < p - 1$, we must have $\theta^\alpha_{1, 1} = \theta^\beta_{0, n} = \theta^\beta_{x, m}$ for some $x \neq 0, m \neq 0$.  Replacing $(t, n)$ by $(x, m)$ if necessary, we may assume that $t \in \bFsp$ in the original equation $\theta^\alpha_{1, 1} = \theta^\beta_{t, n}$.  

Now reindex the sum on the right side of \eqref{eq:eigenvalue_1_1_assumption} by writing $t^{-1}s$ for $s$:
\begin{equation*}
\theta^\beta_{t, n} = \sum_{s \in \bFsp} \zeta_p^{ts}\zeta_{\ell}^{n\beta(s)} = \sum_{s \in \bFsp} \zeta_p^{s}\zeta_{\ell}^{n\beta(t^{-1}s)} = \theta^{n \beta^t}_{1, 1}.
\end{equation*}
Apply the automorphism $\zeta_p \mapsto \zeta_p^a$, as in \eqref{eq:eigenvalue_1_1_galois_act}, to get $\theta^\alpha_{a, 1} = \theta^{n \beta^t}_{a, 1}$ for $a \neq 0$, and then sum over $a \in \bF_p^\times$ to get $\theta^\alpha_{0, 1} = \theta^{n \beta^t}_{0, 1}$ as well.  The equality $\alpha = n \beta^t$ of voltage assignments follows from Lemma \ref{lem:fourier_inverse} and so by Proposition \ref{prop:iso_covers},  $X^{\alpha}$ and $X^{n\beta}$ are isomorphic as $\bZ/\ell$-covering graphs. 
The function $(x, i) \mapsto (t^{-1}x, n^{-1}i)$ provides an isomorphism of graphs $X^\alpha \to X^\beta$.
\end{proof}

\section{The Cayley isomorphism property for $X^\alpha$}\label{sec:CI_property}

In this section, we prove the following rigidity theorem for isomorphisms of $\bZ/\ell$-covers of Paley graphs.  The result, an expanded version of Theorem \ref{thm:rigidity_intro} in the introduction, will be used in the discussion of counterexamples in \S\ref{sec:counterexamples}.
\begin{theorem}\label{thm:CI}
Assume that $q = p^r$ satisfies $r \leq 4$ or $\ell > (q - 1)/2$.
Suppose that the translation invariant $\bZ/\ell$-covers $X^\alpha$ and $X^\beta$ of $X(\bFq)$ are isomorphic as graphs.  Then there exists an isomorphism $f \cn X^\alpha \to X^\beta$ of the form $f(x, i) = (tx^\sigma, ni)$ for some $t \in \bFsquares, \sigma \in \Gal(\bF_q/\bF_p), n \in (\bZ/\ell)^\times$.  Moreover, $n \alpha = \beta^{t^{-1}} \circ \sigma$, and $X^{n \alpha}$ and $X^\beta$ are isomorphic as $\bZ/\ell$-covering graphs.
\end{theorem}

The strategy is to establish a limited form of the Cayley isomorphism property for the group $\bF^+_q \times \bZ/\ell$.  Let $G$ be a finite group.  Consider the Cayley graph $\Cay(G, S)$ of $G$ with respect to a generating set $S \subset G$ closed under taking inverses.  We say that the graph $\Cay(G, S)$ has the Cayley isomorphism property, or more briefly that $S$ is a CI-subset of $G$, if, whenever $\Cay(G, S) \cong \Cay(G, T)$ as graphs, there is a group automorphism $f \cn G \to G$ such that $f(S) = T$.  If all generating sets $S \subset G$ are CI-subsets, then the group $G$ is said to have the CI-property.  Groups known to have the CI-property include $G = \bZ/p$ \citelist{\cite{Elspas_Turner} \cite{Babai:CI}} and more generally $(\bZ/p)^r$ for $r \leq 5$, but not $\bZ/p^r$ when $r > 1$ \citelist{\cite{Babai_Frankl_I} \cite{Godsil} \cite{Feng_Kovacs} \cite{Li_survey}}.  

Determining when $G = \bFplus \times \bZ/\ell \cong (\bZ/p)^r \oplus \bZ/\ell$ has the CI-property is currently an active area \citelist{\cite{Somlai_Muzychuk} \cite{Kovacs_Muzychuk}}.  The best available result, to our knowledge, establishes the CI-property when $r \leq 4$ for any distinct primes $p, \ell$ \cite{Kovacs_Ryabov}, and the hypotheses of Theorem \ref{thm:CI} reflect that status.  In order to cover the cases when the CI-property does not hold for all generating subsets $S$, we take on the modest goal of showing that the generating set $S^\alpha$ associated to a voltage assignment is a CI-subset. More precisely, we prove:

\begin{proposition}\label{prop:CI}
Assume that $\ell > (q - 1)/2$, and let $\alpha \cn \bFsq \to \bZ/\ell$ be a voltage assignment. 
Then the generating set 
\begin{equation*}
S^\alpha = \{(s, \alpha(s)) \mid s \in \bFsq \} 
\end{equation*}
is a CI-subset of $\bF_q^+ \times \bZ/\ell$.
\end{proposition}

The theorem follows immediately from the proposition.  To see this, assume that $X^\alpha$ and $X^\beta$ are isomorphic as graphs and use the Cayley isomorphism property to find a group automorphism $f \cn \bF_q^+ \times \bZ/\ell \to \bF_q^+ \times \bZ/\ell$ satisfying $f(S^\alpha) = S^\beta$.  Such an automorphism takes the form $f(x, i) = (tx^\sigma, ni)$ for some $t \in \bFtimes, \sigma \in \Gal(\bF_q/\bF_p), n \in (\bZ/\ell)^\times$, where we apply Carlitz's theorem \cite{Carlitz} (in its original form) to obtain the description on the first component.  In order that $f(s, \alpha(s)) \in S^\beta$ we must have $t \in \bFsq$ and $\beta(t s^\sigma) = n \alpha(s)$.  It follows that $f$ defines an isomorphism of $\bZ/\ell$-covers $X^{n\alpha} \to X^\beta$ by the description in Proposition \ref{prop:iso_covers}.

\medskip

Before proving the proposition, we first show that the Paley graph $X(\bF_q) = \Cay(\bFplus, \bFsq)$ has the Cayley isomorphism property.  
The proof uses a criterion due to Babai \cite{Babai:CI}*{Lem. 3.1} (see also \cite{Li_survey}*{\S4}): $S$ is a CI-subset of $G$ if and only if every subgroup $H$ of $\Aut(\Cay(G, S))$ that is abstractly isomorphic to $G$ and acts freely and transitively on the set of vertices is conjugate to the subgroup $G < \Aut(\Cay(G, S))$ of translations by elements of the group $G$.  Thus, to establish that $S = \bFsq$ is a CI-subset of $G = \bFplus$, it suffices to prove the following lemma, which is another direct consequence of Carlitz's theorem.

\begin{lemma}\label{lem:Paley_is_CI}
The group $G = \{x \mapsto x + a\}$ of translations is the unique subgroup of $\Aut(X(\bFq))$ isomorphic to $\bF_q^+$.
\end{lemma}

\begin{proof}
Suppose that $H < \Aut(X(\bFq))$ is isomorphic to $\bF_q^+$.  By \eqref{eq:Paley_autos}, the automorphism group is a semidirect product $\Aut(X(\bFq)) = G \ltimes A_0$, where $A_0 = \{ x \mapsto t x^{\sigma}\}$ is the subgroup fixing the vertex $0$.  Since $p$ does not divide $\abs{\bFsq} = (q - 1)/2$, the image of $H \cong (\bZ/p)^r$ under the projection onto $A_0$ must lie in the cyclic subgroup $\Gal(\bF_q/\bF_p) < A_0$ of order $r$.  If the image is non-trivial, then there is an element $\sigma \in \Gal(\bF_q/\bF_p) \cap H$ of order $p$ that commutes with all translations in $G \cap H$.  Such translations are by elements of the fixed field $\bFq^\sigma = \bF_{p^{r/p}}$, and so $p^{r - 1} = \abs{G \cap H} \leq p^{r/p}$, a contradiction.  It follows that $H = G$.
\end{proof}

\begin{proof}[Proof of Proposition \ref{prop:CI}]
We apply the criterion of Babai with $G = \bF_q^+ \times \bZ/\ell$ and $S = S^\alpha$.  Suppose that $H$ is a subgroup of $A = \Aut(X^\alpha)$ isomorphic to $G$ that acts freely and transitively on the vertices of $X^\alpha$.  Let $L$ denote the subgroup $\bZ/\ell < G < A$.  

Consider the case when $L$ is not a Sylow $\ell$-subgroup of $A$.  Then $\ell^2$ divides $\abs{A} = \abs{G} \cdot \abs{A_0}$, where $A_0$ is the stabilizer of a chosen vertex of $X^\alpha$.  Thus, there is an element $f \in A_0$ of order $\ell$.  By the connectivity of $X^\alpha$, there must be an edge 
connecting $x$ and $y$ in $X^\alpha$ with $f(x) = x$ and $f(y) \neq y$.  Then the orbit of $y$ under $f, f^2, \dotsc $ contains $\ell$ vertices, contradicting the fact that the valence of $x$ is $(q - 1)/2 < \ell$.

Now suppose that $L$ is a Sylow $\ell$-subgroup of $A$.  Since $H \cong G$ also contains a Sylow $\ell$-subgroup of $A$, there exists $f \in A$ such that $G \cap fHf^{-1}$ contains $L$.  Note that $G$ and $fHf^{-1}$ are both contained in the centralizer $C_A(L)$ of $L$ in $A$.  The quotient group $C_A(L)/L$ acts on the Paley graph $X(\bF_q)$ under the canonical identification $X^\alpha /L \cong X(\bF_q)$.  The induced action of $G/L$ is the usual translation action of $\bF_q^+$ on $X(\bF_q)$.  The subgroup $fHf^{-1}/L < \Aut(X(\bFq))$ is isomorphic to $\bF_q^+$, so by Lemma \ref{lem:Paley_is_CI} they must coincide: $G/L = fHf^{-1}/L$.  It follows that $G = fHf^{-1}$ as subgroups of $A$.  By Babai's criterion, $S^\alpha$ is a CI-subset of $G = \bF_q^+ \times \bZ/\ell$.  
\end{proof}

\section{The construction of cospectral non-isomorphic covering graphs}\label{sec:counterexamples}

In this section, we construct examples of $\bZ/\ell$-covering graphs $X^\alpha$ and $X^\beta$ of the Paley graph $X(\bF_q)$ that are not isomorphic as graphs yet have the same adjacency spectrum.  Using the notation from \eqref{eq:eigenvalue_formula}, the agreement between the eigenvalues of $X^\alpha$ and $X^\beta$
\begin{equation}\label{eq:goal_spec_bijection}
\theta^{\alpha}_{a, k} = \theta^\beta_{f(a), k} \quad \text{for all $a \in \bF_q$, $k \in \bZ/\ell$,}
\end{equation} 
will be encoded by an odd permutation polynomial $f(T) \in \bF_q[T]$, meaning a polynomial of odd degree that induces a bijection $a \mapsto f(a)$ on the set $\bF_q$.

Given a character $\psi \in (\bZ/\ell)^\vee$ and a voltage assignment $\alpha$, recall the function $\psi^\alpha \cn \bFq \to \bC$ from \S\ref{sec:covers}, which takes the value $\psi(\alpha(s))$ for $s \in \bFsq$ and is otherwise zero.
Next, define an $\bF_q \times \bF_q$-indexed matrix $A(\psi^{\alpha})$ by $A(\psi^{\alpha})_{x, y} = \psi^{\alpha}(y - x)$.
When $\psi = \psi_0$ is the trivial character, the matrix $A(\psi^{\alpha}_0)$ is the adjacency matrix for the Paley graph $X(\bFq)$.  We call $A(\psi^{\alpha})$ the $\psi$-twisted adjacency matrix of $X(\bFq)$ associated to the voltage assignment $\alpha$.  Our standing assumptions on $\alpha$ imply that $A(\psi^{\alpha})$ is a Hermitian matrix and thus has real eigenvalues.   As we will soon see, these eigenvalues are precisely the submultiset $\{\theta^{\alpha}_{a, k}\}_{a \in \bFq}$ of the adjacency spectrum of the covering graph $X^{\alpha}$, where $\psi = \psi_k$.

Notice that the entries of $A(\psi^\alpha)$ remain invariant under additive translations: 
\begin{equation*}
A(\psi^\alpha)_{x + a, y + a} = A(\psi^\alpha)_{x, y} \quad \text{for each $a \in \bFq$.}
\end{equation*}
We call a matrix with this property an $\bFq$-circulant matrix, since in the case $q = p$ it coincides with the well-studied notion of a $p \times p$ circulant matrix \cite{Kra_Simanca}.  The most basic example of an $\bFq$-circulant matrix is the cycle matrix $R^t$ associated to $t \in \bFq$: 
\begin{equation*}
R^t_{x, y} = \begin{cases} 1 \quad &\text{if $y - x = t$,} \\ 0 \quad &\text{else.} \end{cases}
\end{equation*}
Any $\bF_q$-circulant matrix is a linear combination of the matrices $R^t$, and thus is determined by a single $\bF_q$-indexed vector.  For example, the $\psi$-twisted adjacency matrix is given by $A(\psi^\alpha) = \sum_t \psi^\alpha(t)R^t$. Here, and in the arguments to follow, the sum runs over all values $t \in \bFq$, but of course the coefficients $\psi^\alpha(t)$ are only nonzero when $t \in \bFsq$.


\begin{proposition}\label{prop:TFAE}
Suppose given voltage assignments $\alpha$ and $\beta$ encoding a pair of $\bZ/\ell$-covers $X^\alpha$ and $X^\beta$ of $X(\bFq)$.  Fix $k \in \bZ/\ell$ with corresponding character $\psi = \psi_k$. Then the following conditions are equivalent.
\begin{enumerate}  
\item The $\psi$-twisted adjacency matrices $A(\psi^{\alpha})$ and $A(\psi^{\beta})$ are cospectral.
\item There is an odd permutation polynomial $f(T) \in \bF_q[T]$ encoding a bijective correspondence $\theta^{\alpha}_{a, k} = \theta^\beta_{f(a), k}$ between the adjacency eigenvalues of $X^{\alpha}$ and $X^{\beta}$ with second index $k$.
\item There exists an $\bFq \times \bFq$-indexed orthogonal matrix $M$ such that $M \psi^{\beta} = \psi^{\alpha}$ and, for each $t \in \bFq$, the matrix $MR^tM^{\top}$ is the $\bFq$-circulant matrix associated to the column of $M$ indexed by $t$:
\[
MR^tM^{\top} = \sum_{a \in \bFq} M_{a, t} R^a.
\]
\end{enumerate}
When these conditions hold, the matrix $M$ and the permutation polynomial $f(T)$ may be chosen so that 
\begin{equation}\label{def_matrix_M}
M_{x, y} = \frac{1}{q}\sum_{a \in \bF_q}  \chi_{f(a)}(y)\overline{\chi_a(x)}.
\end{equation}
\end{proposition}
\begin{proof}
Suppose given a matrix $M$ satisfying the conditions in (3).  
To prove (1), use the decomposition of the $\bFq$-circulant matrix $A(\psi^\beta) = \sum_{t}\psi^\beta(t)R^t$, and directly calculate that the $\psi$-twisted adjacency matrices are similar:
\begin{align*}
M A(\psi^\beta) M^{-1} = \sum_t \psi^\beta(t) MR^t M^{-1} = \sum_{a, t} M_{a, t} \psi^{\beta}(t) R^a = \sum_a \psi^{\alpha}(a) R^a = A(\psi^\alpha).
\end{align*}
It follows that $A(\psi^\alpha)$ and $A(\psi^\beta)$ have the same spectrum.

Next, we show that (1) implies (2).  A calculation similar to \eqref{eq:char_eigenfunction} shows that each character $\chi \in (\bFplus)^\vee$ is an eigenfunction for $A(\psi^\alpha)$ with eigenvalue $\sum_{t} \chi(t)\psi^\alpha(t)$.  For $\chi = \chi_a$ and $\psi = \psi_k$, this is the eigenvalue $\theta^\alpha_{a, k}$.  Thus, condition (1) is equivalent to the existence of a correspondence of eigenvalues $\theta^{\alpha}_{a, k} = \theta^\beta_{f(a), k}$ for some bijection $f \colon \bFq \to \bFq$, which by a counting argument we may assume to be a polynomial (of degree less than $q$).  The relation $\theta^{\alpha}_{-a, -k} = \theta^{\alpha}_{a, k}$, which follows from \eqref{eq:eigenvalue_formula}, implies that 
\[
\theta^{\alpha}_{-a, k} = \theta^{\alpha}_{a, -k} = \theta^{\beta}_{f(a), -k} = \theta^{\beta}_{-f(a), k}.
\]
Thus $f(-a) = -f(a)$ for all $a \in \bFq$ so $f(T)$ is an odd permutation polynomial.    

To prove that (2) implies (3), let $f(T)$ be an odd permutation polynomial implementing the bijection between the eigenvalues of $X^\alpha$ and $X^\beta$.  Define a matrix $M$ indexed on $\bFq \times \bFq$ by the formula \eqref{def_matrix_M} in the statement of the proposition.
Fix $t \in \bFq$ and consider an arbitrary entry of the matrix $MR^tM^{\top}$:
\[
(MR^tM^{\top})_{x, y} = \sum_{a} M_{x, a}M_{y, a + t} = \frac{1}{q^2} \sum_{a, b, c}\chi_{f(b)}(a) \overline{\chi_{b}(x)} \chi_{f(c)}(a + t)\overline{\chi_{c}(y)} .
\]
Replace $b$ by $-b$ in the sum and use the property $\chi_{f(-b)}(a) = \overline{\chi_{f(b)}(a)}$, which follows from the assumption that $f(T)$ has odd degree, to rewrite the matrix entry as
\begin{align*}
\frac{1}{q^2} \sum_{b, c}  \chi_{b}(x) \chi_{f(c)}(t)\overline{\chi_{c}(y)} \sum_a \chi_{f(c)}(a)\overline{\chi_{f(b)}(a)} &= \frac{1}{q} \sum_{b, c}  \chi_{b}(x) \chi_{f(c)}(t)\overline{\chi_{c}(y)} \langle \chi_{f(c)}, \chi_{f(b)} \rangle.
\end{align*}
By the orthonormality of characters and the bijectivity of $f$, we see that the sum vanishes except for the $b = c$ term, in which case it equals $M_{y - x, t}$.
Thus, $MR^tM^{\top}$ is an $\bFq$-circulant matrix as described in condition (3).  In particular, since the $t = 0$ column $M_{a, 0}$ is all zero except for the entry $M_{0, 0} = 1$, it follows that $MM^{\top} = I$.  Thus, $M$ is an orthogonal matrix and it remains to verify that $M \psi^{\beta} = \psi^\alpha$.  This follows from the identification of eigenvalues $\theta^{\alpha}_{a, k} = \theta^{\beta}_{f(a), k}$ for $\psi = \psi_k$ and another invocation of the orthonormality of characters:
\begin{align*}
\sum_t M_{x, t} \psi^{\beta}(t) &= \frac{1}{q} \sum_{a, t} \chi_{f(a)}(t) \overline{\chi_{a}(x)} \psi^{\beta}(t) \\
&= \frac{1}{q} \sum_{a, t} \chi_{a}(t) \overline{\chi_{a}(x)} \psi^{\alpha}(t) = \psi^{\alpha}(x). \qedhere
\end{align*}
\end{proof}

In practice, the third condition in the proposition can be verified without too much computational effort.  Notice that when the polynomial $f(T)$ can be chosen independently of $k \in \bZ/\ell$, then the matrix $M$ also does not depend on $k$.  In this situation, the proposition gives a bijection between the multisets of adjacency eigenvalues of $X^\alpha$ and $X^\beta$, as in \eqref{eq:goal_spec_bijection}.

\medskip

Before constructing examples, we give a streamlined method for checking (3) when the exponent $r$ in $q = p^r$ is even.  This condition guarantees that elements of $\bF_p^{\times}$ are squares in $\bF_q$, and so $\bFsq$ is partitioned into $\bF_p^{\times}$-orbits.  Suppose that $f(T)$ is an odd permutation polynomial satisfying $f(a T) = a f(T)$ for all $a \in \bF_p$.  Then the associated $(\bFq \times \bFq)$-indexed matrix $M$ defined by \eqref{def_matrix_M} satisfies $M_{ax, a y} = M_{x, y}$ for all $a \in \bF_p^{\times}$.  In other words, $M$ is invariant under multiplicative translation by nonzero scalars in the underlying prime field.  After removing the $x = 0$ row and the $y = 0$ column, any such matrix can be written as a block circulant matrix whose blocks record the interaction between representatives for each $\bF_p^{\times}$-orbit in $\bFq^{\times}$.  To verify the condition $M \psi^{\beta} = \psi^{\alpha}$, we only need to consider the entries $M_{x, y}$ for $(x, y) \in \bFsq \times \bFsq$, so let us write $\resM$ for that submatrix.  


Let $g$ be a multiplicative generator of $\bFtimes$, and order $\bFsq$ by the even powers of $g$.  Reindex the matrix $\resM$ with respect to this ordering, so that $\resM_{i, j} = M_{g^{2i}, g^{2j}}$ for $0 \leq i, j < (q - 1)/2$. 
Observe that $g^{2i} \in \bF_p$ exactly when $i$ is an integer multiple of $\tfrac{q - 1}{2(p - 1)}$.  It follows that $\resM$ is block circulant of the form
\begin{equation}\label{eq:blocks}
{\footnotesize
\resM = \begin{bmatrix} M^{[1]} & M^{[2]} & \dotsm & M^{[p - 1]} \\
	M^{[p - 1]} & M^{[1]} & \dotsm  & M^{[p - 2]} \\
	\vdots & \vdots & \ddots  & \vdots \phantom{\displaystyle \sum} \\
	M^{[2]} & M^{[3]} & \dotsm & M^{[1]} 
\end{bmatrix}
}
\end{equation}
with each block $M^{[k]}$ a $\tfrac{q - 1}{2(p - 1)} \times \tfrac{q - 1}{2(p - 1)}$ matrix.  Using the notation $m(j, k) = \tfrac{(k - 1)(q - 1)}{2(p - 1)} + j$, the blocks are given by $M^{[k]}_{i, j} = \resM_{i, m(j, k)}$. 
The same ordering on $\bFsq$ by even powers of $g$ decomposes $\psi^\alpha$, thought of as an $\bFsq$-indexed column vector, into $p - 1$ blocks $\psi^{\alpha}_{[j]}$ with $\psi^{\alpha}_{[j]}(i) = \psi(\alpha(g^{2m(i, j)}))$ for $0 \leq i < \tfrac{q - 1}{2(p - 1)}$.  It is useful to let the notation $[j]$ for blocks take values modulo $p - 1$.  Then the assumption $\alpha(-s) = - \alpha(s)$ and the fact that $g^{(q - 1)/2} = -1$ imply that $\psi^{\alpha}_{[(p - 1)/2 + j]} = \overline{\psi^{\alpha}_{[j]}}$.  Consequently, when checking the condition $M \psi^{\beta} = \psi^{\alpha}$ in (3) of Prop. \ref{prop:TFAE}, it suffices to verify that 
\begin{equation}\label{eq:to_check}
\sum_{i = 1}^{p - 1} M^{[i]} \psi^{\beta}_{[i + j - 1]} = \psi^{\alpha}_{[j]} \qquad \text{for $j = 1, \dotsc, (p - 1)/2$,}
\end{equation}
since the remaining equations for the other $(p - 1)/2$ blocks of $\psi^{\alpha}$ follow by applying complex conjugation.

\begin{example}\label{ex:main_counterexample}
Let $q = 25$, and consider the Paley graph $X(\bF_{25})$.  We will construct non-isomorphic  $\bZ/\ell$-covers $X^{\alpha}, X^{\beta}$ of $X(\bF_{25})$ with the same adjacency spectrum.  The method works for any odd prime $\ell$.

Write $\bF_{25} = \bF_{5}(\omega)$ with $\omega^2 = 2$, and observe that the squares in $\bF_{25}$ consist of the $\bF_5^\times$-orbits of $1, \omega + 1$, and $4\omega + 1$.  
Define voltage assignments $\alpha, \beta \cn \bF_{25}^{\boxtimes} \to \bZ/\ell$ by 
\begin{equation*}
\begin{array}{ r|cccccc }
	s  & 1 & 3 & \omega + 1 & 3 \omega + 3 & 4\omega + 1 & 2 \omega + 3 \\ \hline
	\alpha(s) & 1 & 1 & 1 & 0 & 0 & 0 \\
	\beta(s) & 0 & 0 & 1 & 0 & 1 & 1 
\end{array}
\end{equation*}
where we use the relations $\alpha(-s) = -\alpha(s)$ and $\beta(-s) = - \beta(s)$ to determine the remaining values.  In other words, we first define the voltage assignments on $s \in \bF_{25}^{\boxtimes}$ with $\tr(s) = 1, 2$ and then the extension to squares of trace class $3, 4$ is uniquely determined.  Also note that $\alpha$ and $\beta$ are related by $\alpha(s) = \beta(\tau(s))$, where $\tau$ is the trace-preserving permutation of $\bF_{25}^{\boxtimes}$ that swaps the orbits of $1$ and $4\omega + 1$ while fixing the orbit of $\omega + 1$.  This is similar to, but different from, the non-trivial Galois automorphism $\omega \mapsto -\omega$, which also preserves trace, but swaps the orbits of $\omega + 1$ and $4\omega + 1$ while fixing the orbit of $1$.  

We first prove that there is no isomorphism of graphs $X^{\alpha} \cong X^{\beta}$.  If there were, then by Theorem \ref{thm:CI} we would have $n\alpha(s) = \beta(t s^\sigma)$ for some $n \in (\bZ/\ell)^\times$, $t \in \bF_{25}^\times$, $\sigma \in \Gal(\bF_{25}/\bF_{5})$.  It follows that $n = n \alpha(1) = \beta(t)$, so we must have $n = \pm 1$ by the possible values of $\beta(t)$.  This gives $\alpha(s) = \beta(\pm t s^{\sigma})$.  We cannot directly conclude that $\tau(s) = \pm t s^{\sigma}$, which would be impossible because of the structure of $\tau$.  But a case by case argument is still accessible.  Since $\beta(\pm 1) = 1$, the only possibilities are $\pm t = \omega + 1, 4\omega + 1, 2\omega + 3$.  When $\pm t = \omega + 1$, then $s = 3$ gives $1 = \alpha(3) = \beta(3\omega + 3) = 0$, a contradiction.  Similar calculations with $s = \omega + 1$ exclude all other possible values of $\pm t$, including when $\sigma$ is both trivial and non-trivial.  

Having shown that $X^{\alpha}$ and $X^{\beta}$ are not isomorphic as graphs, we next verify that their adjacency spectra are equivalent.  To this end, consider the polynomial
\[
f(T) = (3\omega + 4)T^{21} + (\omega + 1)T^{17} + (4\omega + 1)T^{13} + (2\omega + 4)T^9 + (4\omega + 4)T^5 + (\omega + 2)T.
\]
It can be directly verified using a computer that $f(T)$ is an odd permutation polynomial encoding a bijection between the adjacency eigenvalues of $X^\alpha$ and $X^\beta$, but we prefer to use properties of $f(T)$ to check condition (3) in Proposition \ref{prop:TFAE} for the $(\bF_{25} \times \bF_{25})$-indexed matrix
\[
M_{x, y} = \frac{1}{25}\sum_{a \in \bF_{25}}  \chi_{f(a)}(y)\overline{\chi_a(x)} = \frac{1}{25} \sum_{a \in \bF_{25}}  \zeta_5^{\tr(f(a)y - ax)}.
\]
As in the proof of the proposition, it is immediate from the formula that $M$ is an orthogonal matrix that is $\bF_q$-circulant of the form given in condition (3), so it remains only to check that $M\psi^\alpha = \psi^\beta$ for each character $\psi$ of $\bZ/\ell$.

First note that the coefficients of $f(T)$ sum to one, and so $f(1) = 1$.  The powers of $T$ are arranged so that $f(a T) = a f(T)$ for all $a \in \bF_5$ and so $\resM$ is a block circulant matrix of the form \eqref{eq:blocks}.  In the calculations below, we have used the choice $g = \omega + 3$ for a multiplicative generator.  With respect to the ordering of $\bF_{25}^{\times}$ by even powers of $g$, the blocks of $M$ are	
{\footnotesize
\begin{equation*}
M^{[1]} = \dfrac{1}{5}
\begin{bmatrix}
	0 & 1 & 1 \\
	0 & 3 & 1 \\
	1 & -1 & 1
\end{bmatrix}
\;
M^{[2]} = \dfrac{1}{5}
\begin{bmatrix}
	0 & 1 & 3 \\
	1 & 0 & 1 \\
	2 & 0 & 0
\end{bmatrix}
\;
M^{[3]} = \dfrac{1}{5}
\begin{bmatrix}
	-1 & -1 & 1 \\
	0 & 2 & -1 \\
	1 & 1 & 0
\end{bmatrix}
\;
M^{[4]} = \dfrac{1}{5}
\begin{bmatrix}
	1 & -1 & 0 \\
	-1 & 0 & -1 \\
	1 & 0 & -1
\end{bmatrix}.
\end{equation*}
}Set $\zeta = \psi(1)$.  Then the two cases of equation \eqref{eq:to_check} that must be checked are given by the calculations
{\footnotesize
\begin{equation*}
	\dfrac{1}{5}
\begin{bmatrix}
	0 & 1 & 1 \\
	0 & 3 & 1 \\
	1 & -1 & 1
\end{bmatrix}
\begin{bmatrix}
	1 \\
	\zeta \\
	\zeta
\end{bmatrix}
+
\dfrac{1}{5}
\begin{bmatrix}
	0 & 1 & 3 \\
	1 & 0 & 1 \\
	2 & 0 & 0
\end{bmatrix}
\begin{bmatrix}
	1 \\
	1 \\
	\zeta
\end{bmatrix}
+
\dfrac{1}{5}
\begin{bmatrix}
	-1 & -1 & 1 \\
	0 & 2 & -1 \\
	1 & 1 & 0
\end{bmatrix}
\begin{bmatrix}
	1 \\
	\zeta^{-1} \\
	\zeta^{-1}
\end{bmatrix}
+
\dfrac{1}{5}
\begin{bmatrix}
	1 & -1 & 0 \\
	-1 & 0 & -1 \\
	1 & 0 & -1
\end{bmatrix}
\begin{bmatrix}
	1 \\
	1 \\
	\zeta^{-1}
\end{bmatrix}
=
\begin{bmatrix}
	\zeta \\
	\zeta \\
	1
\end{bmatrix}
\end{equation*}
}

{\footnotesize
\begin{equation*}
	\dfrac{1}{5}
\begin{bmatrix}
	1 & -1 & 0 \\
	-1 & 0 & -1 \\
	1 & 0 & -1
\end{bmatrix}
\begin{bmatrix}
	1 \\
	\zeta \\
	\zeta
\end{bmatrix}
+
\dfrac{1}{5}
\begin{bmatrix}
	0 & 1 & 1 \\
	0 & 3 & 1 \\
	1 & -1 & 1
\end{bmatrix}
\begin{bmatrix}
	1 \\
	1 \\
	\zeta
\end{bmatrix}
+
\dfrac{1}{5}
\begin{bmatrix}
	0 & 1 & 3 \\
	1 & 0 & 1 \\
	2 & 0 & 0
\end{bmatrix}
\begin{bmatrix}
	1 \\
	\zeta^{-1} \\
	\zeta^{-1}
\end{bmatrix}
+
\dfrac{1}{5}
\begin{bmatrix}
	-1 & -1 & 1 \\
	0 & 2 & -1 \\
	1 & 1 & 0
\end{bmatrix}
\begin{bmatrix}
	1 \\
	1 \\
	\zeta^{-1}
\end{bmatrix}
=
\begin{bmatrix}
	\zeta^{-1} \\
	1 \\
	1
\end{bmatrix}.
\end{equation*}
}
Notice that the case of the trivial character $\psi = \psi_0$ is included because the calculations are valid for $\zeta = 1$.  But this case also follows by observing that 
{\footnotesize
\begin{equation*}
\sum_{i = 1}^{4} M^{[i]} = \begin{bmatrix}
	0 & 0 & 1 \\
	0 & 1 & 0 \\
	1 & 0 & 0
\end{bmatrix}
\end{equation*}
}fixes the all-$1$s vector.  It is interesting to observe that the permutation matrix arising here is a manifestation of the permutation $\tau$ of $\bF_{25}^{\boxtimes}$ used in constructing $\alpha$ and $\beta$.
\end{example}

\begin{figure}[h]
\begin{tabular}{ccc}
\includegraphics[scale=.6]{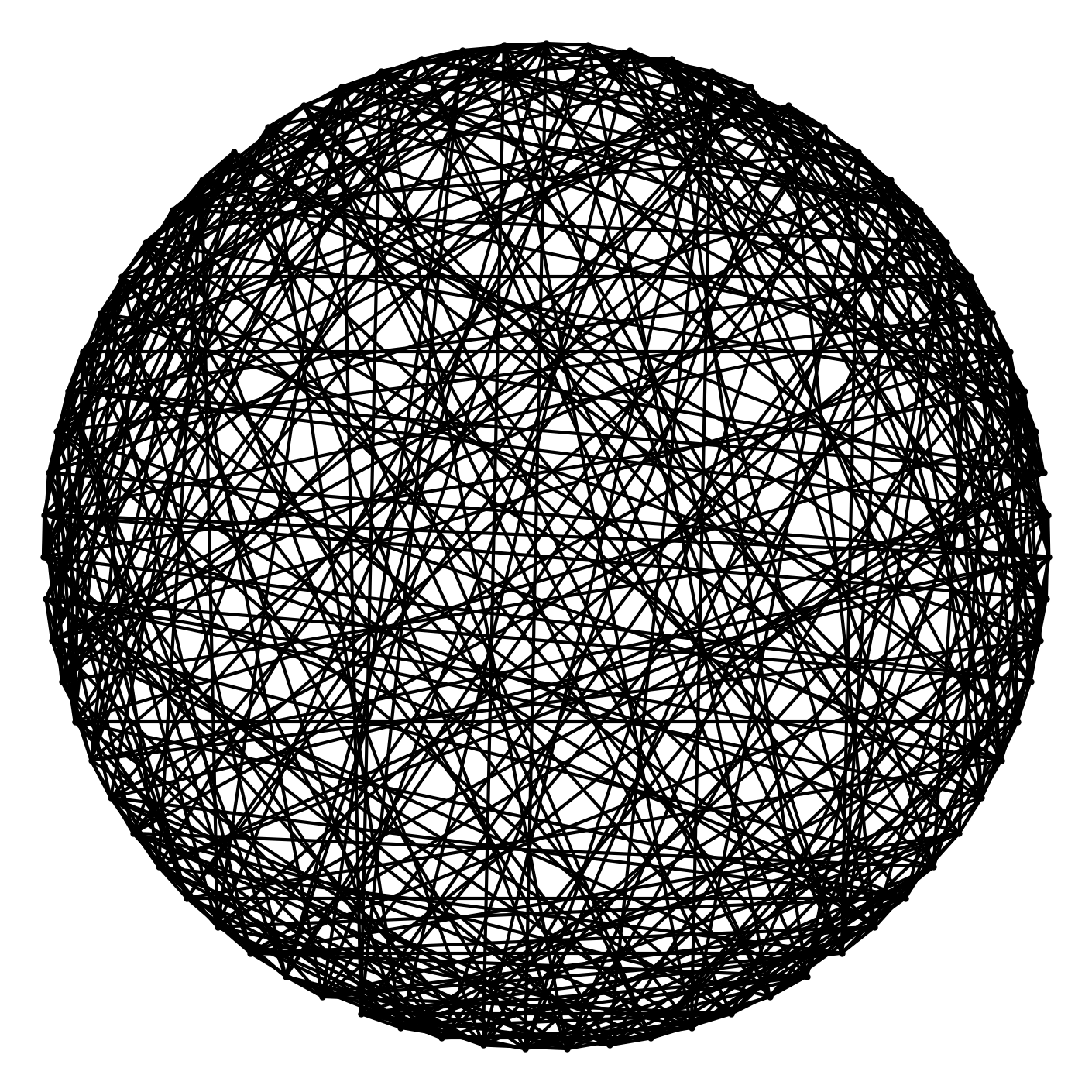} & \qquad &  \includegraphics[scale=.6]{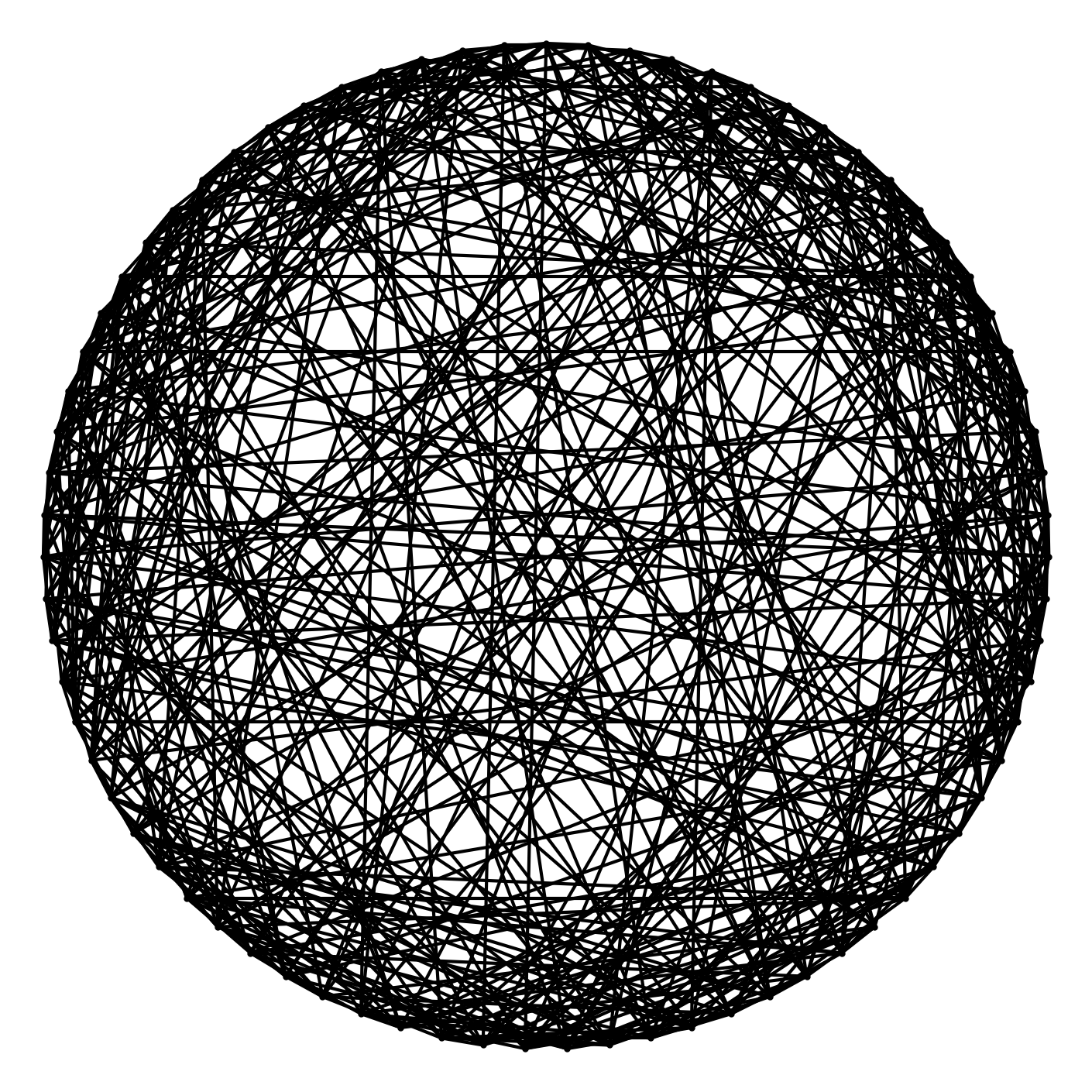} 
\end{tabular}
\caption{A minimal example of cospectral nonisomorphic covers ($q = 25$, $\ell = 3$)}
\end{figure}

A similar technique, using a trace-preserving permutation $\tau$ of order two that breaks Galois symmetry, can be used to construct cospectral non-isomorphic $\bZ/\ell$-covers of the Paley graph $X(\bF_q)$ when $q = p^2$ for any prime $p > 3$.  Computer calculations indicate that such counterexamples exist for all $q = p^r$ with $r > 1$ except $q = 9$.  It would be interesting to parametrize the phenomenon in a systematic way.

\bibliographystyle{plain} 
\bibliography{references.bib}

\end{document}